\documentclass{amsart}
\usepackage{amssymb,graphicx}

\swapnumbers
\newtheorem{theorem}[subsubsection]{Theorem}
\newtheorem{proposition}[subsubsection]{Proposition}
\newtheorem{lemma}[subsubsection]{Lemma}
\theoremstyle{remark}
\newtheorem{remark}[subsubsection]{Remark}
\numberwithin{equation}{section}

\def\ie{\emph{i.e.}}
\def\via{\emph{via}}
\def\cf{\emph{cf}}
\def\eg{\emph{e.g.}}

{\catcode`\@=11
\gdef\mnote#1{\marginpar{\footnotesize
 \tolerance\@M\spaceskip2.6\p@ plus10\p@ minus.9\p@\rm#1}}}
\def\Dg:{\endgraf{\bf Dg:\enspace}\ignorespaces}

\let\bold\mathbf
\let\Bbb\mathbb
\let\Cal\mathcal

\let\sA\bA 
\let\sD\bD 
\let\sE\bE 
\def\DG#1{\Bbb D_{#1}} 
\def\CG#1{\Z_{#1}}     
\let\curve=B  

\let\BB\Curve
\let\sminus\smallsetminus 
\let\splus\oplus 
\let\ge\geqslant 
\let\le\leqslant 
\let\1e 
\def\rel(#1):#2\endrel{$$
 \hbox to\displaywidth{$(#1):$\hfill$\displaystyle#2$\hfill}$$}

\def\PP{\bar P}
\def\QQ{\bar Q}
\def\SS{\bar S}

\def\tX{\tilde X}
\def\tSigma{\tilde\Sigma}
\def\tc{\tilde c}

\def\PGL{\operatorname{\text{\sl PGL}}}
\def\GL{\operatorname{\text{\sl GL}}}

\let\Ge\epsilon
\let\Gg\gamma
\let\Go\omega
\let\eps\varepsilon

\def\Z{\Bbb Z}
\def\R{\Bbb R}
\def\C{\Bbb C}
\def\Q{\Bbb Q}

\def\F{\Bbb F}
\def\CL{\Cal L}
\def\CK{\Cal K}
\def\Cp#1{\Bbb P^{#1}}

\def\Hom{\operatorname{Hom}}
\def\discr{\operatorname{discr}}

\def\id{\operatorname{id}}
\def\Pic{\operatorname{Pic}}
\def\ord{\operatorname{ord}}

\def\Term#1{$\DG{#1}$}
\def\term{\Term{14}}

\def\rom#1{\/{\rm #1}}

\let\PLUS+
{\catcode`\*\active \catcode`\[\active \catcode`\+\active
\let\+\relax
\gdef\deMaple{\catcode`\*\active \catcode`\[\active \catcode`\+\active
\def*{}\def[##1]{_{##1}}\def+{{}\PLUS{}}}}

\begin{document}

\title{A plane sextic with finite fundamental group}

\author[A.~Degtyarev]{Alex Degtyarev}
\address{A.~Degtyarev:
 Department of Mathematics, Bilkent University\\
Bilkent, Ankara 06533, Turkey}
\email{degt@fen.bilkent.edu.tr}

\author[M.~Oka]{Mutsuo Oka}
\address{M.~Oka: Department of Mathematics,
 Tokyo  University of Science\\
 26 Wakamiya-cho, Shinjuku-ku, Tokyo 162-8601, Japan
}
\email{oka@rs.kagu.tus.ac.jp}



\date{}

\dedicatory{}

\begin{abstract}
We analyze irreducible plane sextics whose fundamental group
factors to~$\DG{14}$. We produce explicit equations for all curves
and show that, in the simplest case of the set of singularities
$3\sA_6$, the group is $\DG{14}\times\CG3$.
\end{abstract}

\maketitle

\section{Introduction}

\subsection{Motivation and principal results}
In this paper, we use the term \emph{\Term{2n}-sextic}
for
an irreducible plane sextic $B\subset\Cp2$
whose fundamental group $\pi_1(\Cp2\sminus B)$
factors to the dihedral group~$\DG{2n}$, $n\ge3$.
The \Term{2n}-sextics are classified in~\cite{degt.Oka}
and~\cite{degt.Oka2} (the case of non-simple singularities). The
integer~$n$ can take values~$3$, $5$, and~$7$.
All \Term{6}-sextics are of torus type (\ie, they are given by
equations of the form $p^3+q^2=0$); in particular, their
fundamental groups are infinite.
The \Term{10}-sextics form $13$ equisingular deformation families,
and their fundamental groups are known, see~\cite{degt.D10}
and~\cite{Oka.D10};
with one exception, they are all finite.
Finally, the \term-sextics form two equisingular families, see
Proposition~\ref{prop.sing} below, and their groups are not known.
In this paper, we compute one of the two groups.
Our principal result is the following statement.

\begin{theorem}\label{th.main}
The fundamental group of the complement of
a \term-sextic with the set of
singularities $3\sA_6$ is $\DG{14}\times\CG3$.
\end{theorem}

Theorem~\ref{th.main} is proved in Section~\ref{proof.main}.

Our result can be regarded as another
attempt to substantiate a modified version~\cite{degt.D10} of
Oka's conjecture~\cite{Oka.conjecture} on the fundamental group of
an irreducible plane sextic, stating that the group of an
irreducible sextic with simple singularities that is not of torus
type is finite. (Note that the finiteness of the group is sufficient
to conclude that the Alexander polynomial of the curve is trivial,
see, \eg,~\cite{Oka.survey}.)

\subsection{Contents of the paper}
In Section~\ref{S.construction} we analyze the geometric
properties of \term-sextics, whose existence was proved
in~\cite{degt.Oka} purely arithmetically. We use the theory of
$K3$-surfaces to show that any \term-sextic admits a
$\CG3$-symmetry, see Theorem~\ref{th.symmetry}, and we use this
symmetry to obtain explicit equations defining all \term-sextics,
see Theorem~\ref{th.equation}. The curves form a dimension one
family, depending on one parameter $t\in\C$, $t^3\ne1$.
Most calculations involving
polynomials were done using {\tt Maple}.

The heart of the paper is Section~\ref{S.group}. We use a
particular value $t=5/6$ of the parameter (close to $t=1$, where
the curve degenerates to a triple cubic) and analyze the real part
of the curve obtained. With respect to an appropriately chosen
real pencil of lines, it has sufficiently many real critical
values, and we apply van Kampen's method (ignoring all non-real
critical values)
to produce an `upper
estimate' on the fundamental group, see Theorem~\ref{th.epi}.
Comparing the latter with the known `lower estimate' (the fact
that the curve is known to be a \term-sextic), we prove
Theorem~\ref{th.main}.


\section{The construction}\label{S.construction}

\subsection{Statements}
A $\CG3$-action on~$\Cp2$ is called \emph{regular} if it lifts to
a regular representation $\CG3\to\GL(3,\C)$. An order~$3$ element
$c\in\PGL(3,\C)$ is called \emph{regular} if it generates a
regular $\CG3$-action. Any regular order~$3$ automorphism
of~$\Cp2$ has three isolated fixed points (and no other fixed
points). Conversely, any order~$3$ automorphism~$c$ of~$\Cp2$ with
isolated fixed points only is regular (as isolated fixed points
correspond to dimension one eigenspaces of the lift of~$c$
to~$\C^3$).

The following statement is proved in~\cite{degt.Oka} (see
also~\cite{degt.Oka2}, where sextics with non-simple singular
points are ruled out).

\begin{proposition}\label{prop.sing}
All \term-sextics form two equisingular deformation families, one
family for each of the sets of singularities~$3\sA_6$ and
$3\sA_6\splus\sA_1$.
\qed
\end{proposition}

The principal results of this section are the following two
theorems.

\begin{theorem}\label{th.symmetry}
Any \term-sextic~$B$
is invariant under a certain regular order~$3$ automorphism
$c\colon\Cp2\to\Cp2$ acting
on the three type~$\sA_6$ singular
points of~$B$ by a cyclic permutation.
\end{theorem}

\begin{theorem}\label{th.equation}
Up to projective transformation, the \term-sextics
form a connected one parameter family $\curve(t)$\rom;
in
appropriate homogeneous coordinates, they are given by the
polynomial
\begin{equation}
\deMaple\obeylines\def
{\\}%
\gathered%
2*t*(t^3-1)*(z[0]^4*z[1]*z[2]+z[1]^4*z[2]*z[0]+z[2]^4*z[0]*z[1])
+(t^3-1)*(z[0]^4*z[1]^2+z[1]^4*z[2]^2+z[2]^4*z[0]^2)
+t^2*(t^3-1)*(z[0]^4*z[2]^2+z[1]^4*z[0]^2+z[2]^4*z[1]^2)
+2*t*(t^3+1)*(z[0]^3*z[1]^3+z[1]^3*z[2]^3+z[2]^3*z[0]^3)
+4*t^2*(t^3+2)*(z[0]^3*z[1]^2*z[2]+z[1]^3*z[2]^2*z[0]+z[2]^3*z[0]^2*z[1])
+2*(t^6+4*t^3+1)*(z[0]^3*z[1]*z[2]^2+z[1]^3*z[2]*z[0]^2+z[2]^3*z[0]*z[1]^2)
+t*(t^6+13*t^3+10)*z[0]^2*z[1]^2*z[2]^2, \endgathered\label{equation}%
\end{equation}
where $t\in\C$ and $t^3\ne1$. The restriction
of~$\curve(t)$ to the subset $t^3\ne1,-27$ is an equisingular
deformation, all curves having the set of singularities $3\sA_6$.
The three curves with $t^3=-27$ are extra singular\rom; their sets
of singularities are $3\sA_6\splus\sA_1$.
\end{theorem}

\begin{remark}\label{rem.triple}
We do not assert that all curves~$\curve(t)$
are
pairwise distinct. In fact, one can observe that the substitution
$t\mapsto\Ge t$, $\Ge^3=1$, results in an equivalent
curve, the corresponding change of coordinates being
$(z_0:z_1:z_2)\mapsto(z_0:\Ge^2z_1:\Ge z_2)$. In particular, all
three extra singular curves are equivalent.
\end{remark}

Theorems~\ref{th.symmetry} and~\ref{th.equation} are proved,
respectively,
in Sections~\ref{proof.symmetry} and~\ref{proof.equation} below.

\subsection{Discriminant forms}\label{s.discr}
An \emph{\rom(integral\rom) lattice} is a finitely generated free
abelian group~$L$ supplied with a symmetric bilinear form $b\colon
L\otimes L\to\Z$. We abbreviate $b(x,y)=x\cdot y$ and
$b(x,x)=x^2$. A lattice~$L$ is \emph{even} if $x^2=0\bmod2$ for
all $x\in L$. As the transition matrix between two integral bases
has determinant $\pm1$, the determinant $\det L\in\Z$ (\ie, the
determinant of the Gram matrix of~$b$ in any basis of~$L$) is well
defined. A lattice~$L$ is called \emph{nondegenerate} if $\det
L\ne0$; it is called \emph{unimodular} if $\det L=\pm1$.

Given a lattice~$L$,
the bilinear form extends to $L\otimes\Q$ by linearity. If
$L$ is nondegenerate, the dual group $L^*=\Hom(L,\Z)$ can
be identified with the subgroup
$$
\bigl\{x\in L\otimes\Q\bigm|
 \text{$x\cdot y\in\Z$ for all $x\in L$}\bigr\}.
$$
In particular, $L\subset L^*$. The quotient $L^*/L$ is a finite
group; it is called the \emph{discriminant group} of~$L$ and is
denoted by $\discr L$ or~$\CL$. The discriminant group~$\CL$
inherits from $L\otimes\Q$ a symmetric bilinear form
$\CL\otimes\CL\to\Q/\Z$, called the \emph{discriminant form}, and,
if $L$ is even, its quadratic extension $\CL\to\Q/2\Z$. When
speaking about the discriminant groups, their (anti-)isomorphisms,
etc., we always assume that the discriminant form (and its
quadratic extension if the lattice is even) is taken into account.
One has $\#\CL=\mathopen|\det L\mathclose|$; in particular,
$\CL=0$ if and only if $L$ is unimodular.

An \emph{extension} of a lattice~$L$ is another lattice~$M$
containing~$L$, so that the form on~$L$ is the restriction of that
on~$M$. An \emph{isomorphism} between two extensions
$M_1\supset L$ and $M_2\supset L$ is an isometry $M_1\to M_2$
whose restriction to~$L$ is the identity.
In what follows, we are only interested in the case when
both~$L$ and~$M$ are even and $[M:L]<\infty$. Next
two theorems are found in V.~V.~Nikulin~\cite{Nikulin}.

\begin{theorem}\label{th.Nik1}
Given a nondegenerate even
lattice~$L$, there is a canonical one-to-one correspondence
between the set of isomorphism classes of finite index extensions
$M\supset L$ \rom(by even lattices\rom)
and the set of isotropic subgroups $\CK\subset\CL$.
Under this correspondence,
$M=\{x\in L^*\,|\,x\bmod L\in\CK\}$ and
$\discr M=\CK^\perp\!/\CK$.
\qed
\end{theorem}

The isotropic subgroup $\CK\subset\CL$ as in Theorem~\ref{th.Nik1}
is called the \emph{kernel} of the extension $M\supset L$. It can
be defined as the image of $M\!/L$ under the homomorphism induced by
the natural inclusion $M\hookrightarrow L^*$.

\begin{theorem}\label{th.Nik2}
Let $M\supset L$ be a finite index
extension of a nondegenerate even lattice~$L$
\rom(by an even lattice~$M$\rom), and let $\CK\subset\CL$ be its
kernel.
Then, an auto-isometry $L\to L$ extends to~$M$ if and only if the
induced automorphism of~$\CL$ preserves~$\CK$.
\qed
\end{theorem}

\subsection{Proof of Theorem~\ref{th.symmetry}}\label{proof.symmetry}
Fix a \term-sextic $B\subset\Cp2$ and consider the double covering
$X\to\Cp2$ ramified at~$B$ and its minimal resolution~$\tX$. Since
all singular points of~$B$ are simple, see Proposition~\ref{prop.sing},
$\tX$ is a $K3$-surface.
For a singular point~$P$ of~$B$, denote by~$D_P$ the set of
exceptional divisors in~$\tX$ over~$P$, as well as its incidence
graph, which is the Dynkin graph of the same name
$\sA$--$\sD$--$\sE$ as~$P$.
Let $\Sigma_P\subset H_2(\tX)$ be the sublattice spanned
by~$D_P$. (Here, $H_2(\tX)$ is regarded as a lattice \via\ the
intersection index form.)
Let, further,
$\Sigma'=\bigoplus_P\Sigma_P$, the summation running over all
type~$\sA_6$ singular points~$P$ of~$B$
(see Proposition~\ref{prop.sing}),
and let $\tSigma'\supset\Sigma'$ be the
primitive hull of~$\Sigma'$ in $H_2(\tX)$, \ie,
$\tSigma'=(\Sigma'\otimes\Q)\cap H_2(\tX)$.
It is a finite index
extension; denote by $\CK\subset\discr\Sigma'$ its kernel.

Let $P_0$,
$P_1$, $P_2$ be the type~$\sA_6$ points. For each point $P=P_i$,
$i=0,1,2$, fix an orientation of its (linear) graph~$D_P$ and let
$e_{i1},\ldots,e_{i6}$ be the elements of~$D_P$ numbered
consecutively according to the chosen orientation. Denote by
$e^*_{i1},\ldots,e^*_{i6}$ the dual basis for $\Sigma^*_P$. The
discriminant group $\discr\Sigma_P\cong\Z_7$ is generated by
$e^*_{i1}\bmod\Sigma_P$, and, for each $k=1,\ldots,6$, one
has $e^*_{ik}=ke^*_{i1}\bmod\Sigma_P$. Let
\begin{equation}
\Gg_0=e^*_{04}+e^*_{12}+e^*_{21},\quad
\Gg_1=e^*_{01}+e^*_{14}+e^*_{22},\quad
\Gg_2=e^*_{02}+e^*_{11}+e^*_{24}.\label{eq.gamma}
\end{equation}
According to~\cite{degt.Oka}, under an appropriate numbering
of the type~$\sA_6$ singular points of~$B$ and appropriate
orientation of their Dynkin graphs,
the kernel $\CK\cong\Z_7$ is
generated by the residue $\Gg_0\bmod\Sigma'$. (For the convenience
of the further exposition, we use an indexing slightly different
from that used in~\cite{degt.Oka}.) Observe that each of the
residues
$\Gg_1=2\Gg_0\bmod\Sigma'$ and $\Gg_2=4\Gg_0\bmod\Sigma'$ also
generates~$\CK$.

Define an isometry $c_\Sigma\colon\Sigma'\to\Sigma'$ \via\
$e_{0k}\mapsto e_{1k}\mapsto e_{2k}\mapsto e_{0k}$, $k=1,\ldots6$.
Clearly, $c_\Sigma^3=\id$ and the induced action on
$\discr\Sigma'$ is a regular representation of~$\CG3$ over~$\F_7$.
Hence, $\discr\Sigma'$ splits into direct (not orthogonal) sum of
1-dimensional eigenspaces, $\discr\Sigma'=V_1\oplus V_2\oplus V_4$,
corresponding to the three cubic roots of unity
$1,2,4\in\F_7$.
Since $c_\Sigma(\Gg_0)=\Gg_1=2\Gg_0\bmod\Sigma'$, see above,
one has $\CK=V_2$
and it is immediate that $\CK^\perp\!/\CK$ can be
identified with~$V_1$. Hence, $c_\Sigma$ extends to an
auto-isometry $\tc_\Sigma\colon\tSigma'\to\tSigma'$, see
Theorem~\ref{th.Nik2},
and the induced action on $\discr\tSigma=V_1$, see
Theorem~\ref{th.Nik1}, is trivial. Applying Theorem~\ref{th.Nik2}
to the finite index extension
$H_2(\tX)\supset\tSigma'\oplus(\tSigma')^\perp$, one concludes
that the direct sum $\tc_\Sigma\oplus\id$ extends to an order~$3$
auto-isometry $\tc_*\colon H_2(\tX)\to H_2(\tX)$.

By construction, $\tc_*$ preserves the class~$h$ of the pull-back
of a generic line in~$\Cp2$ and the class~$\Go$ of a holomorphic
$2$-form on~$\tX$ (as both $h,\Go\in(\tSigma')^\perp$).
Furthermore, $\tc_*$ preserve the positive cone~$V^+$ of~$\tX$.
(Recall that the positive cone is an open
fundamental
polyhedron $V^+\subset(\Pic X)\otimes\R$ of the group generated by
reflections defined by vectors $x\in\Pic X$ with $x^2=-2$; in the
case under consideration,
it is uniquely
characterized by the requirement that $V^+\cdot e>0$ for any
exceptional divisor~$e$ over a singular point of~$B$
and that the closure of~$V^+$ should contain~$h$.) The usual
averaging argument shows that $\tX$ has a K\"ahler metric with
$\tc_*$-invariant fundamental class $\rho\in V^+$. The pair
$(\Go\bmod\C^*,\rho\bmod\R^*)$ represents a point in the fine
period space of marked quasipolarized $K3$-surfaces, see
A.~Beauville~\cite{Beauville}, and, since this point is fixed
by~$\tc_*$, there is a unique automorphism $\tc\colon\tX\to\tX$
inducing~$\tc_*$ in the homology. It is of order~$3$ (as
the only automorphism inducing $\tc_*^3=\id$ is the identity),
symplectic (\ie, preserving holomorphic $2$-forms), and commutes
with the deck translation of the ramified covering $\tX\to\Cp2$
(as the map $\tX\to\Cp2$ is defined by the linear system
$h\in\Pic\tX$ preserved by~$\tc_*$).
Thus, $\tc$ descends to an
order~$3$ automorphism $c\colon\Cp2\to\Cp2$. The latter
preserves~$B$ (as it lifts to~$\tX$) and has isolated
fixed points only (as so does its lift~$\tc$, as any symplectic
automorphism of a $K3$-surface); in particular, $c$ is regular.
\qed

\subsection{Geometric properties of \term-sextics}\label{s.basic}
The following geometric characterization of \term-sextics is found
in~\cite{degt.Oka}.

\begin{proposition}\label{prop.conics}
The three type~$\sA_6$ singular points~$P_0$, $P_1$, $P_2$ of a
\term-sextic~$B$ can be ordered so that there are
three conics
$Q_0$, $Q_1$, $Q_2$ such that each $Q_i$,
$i=0,1,2$,
intersects~$B$ at~$P_{i-k}$, $k=1,2,3$, with multiplicity~$2k$.
\end{proposition}

\begin{remark}
Here and below, to shorten the notation, we use the cyclic indexing
$P_{i+3s}=P_i$ and $Q_{i+3s}=Q_i$ for $s\in\Z$.
In fact, the points should be ordered as explained in
Section~\ref{proof.symmetry}; then $Q_i$ is the projection
to~$\Cp2$ of the rational curve realizing the $(-2)$-class
$\Gg_i+h\in\Pic\tX$, where $\Gg_i$ is given by~\eqref{eq.gamma}.
\end{remark}

\begin{lemma}\label{conics}
The automorphism~$c$ given by Theorem~\ref{th.symmetry} acts
on the set of conics $\{Q_0,Q_1,Q_2\}$
as in Proposition~\ref{prop.conics} by a cyclic permutation.
\end{lemma}

\begin{proof}
For each $i=0,1,2$, the incidence conditions described above
define at most one conic~$Q_i$ (as otherwise two conics would
intersect at six points). Since $c$ permutes the singular points
of~$B$, it must also permute the conics.
\end{proof}

\begin{lemma}\label{irreducible}
Let $Q_0$, $Q_1$, $Q_2$ be the conics as
in Proposition~\ref{prop.conics}.
Then, either all~$Q_i$, $i=0,1,2$, are irreducible
or else,
for each $i=0,1,2$, one has a splitting
$Q_i=(P_{i-1}P_i)+(P_iP_{i+1})$. In the latter case, $B$ is
tangent to $(P_{i-1}P_i)$ at~$P_i$.
\end{lemma}

\begin{proof}
Due to Lemma~\ref{conics}, if one
of~$Q_i$ is reducible, so are the others.
Assume that $Q_0$ splits into two lines, $Q_0=L_0'+L_0''$. If the
intersection point $L_0'\cap L_0''$ is a singular point of~$B$,
one immediately concludes that
$Q_0=(P_2P_0)+(P_0P_1)$ and extends this splitting to the other
conics \via~$c$.

Otherwise,
assume that it is~$L_0'$ that
intersects~$B$ at~$P_0$ with multiplicity~$6$. Then the component
$L_2''=c^2(L_0'')$ of $Q_2=c^2(Q_0)$ is tangent to~$B$ at~$P_0$
(we assume that $c$ acts \via\
$P_0\mapsto P_1\mapsto P_2\mapsto P_0$);
hence, $L_2''=L_0'$ and this line cannot pass through~$P_1$.
(Neither can the other component $L_2'=c^2(L_0')$, as it
intersects~$B$ at~$P_2$ with the maximal multiplicity~$6$.)
\end{proof}

\subsection{Theorem~\ref{th.equation}: the generic case}\label{s.generic}
Fix a \term-sextic~$B$ and denote by $P_0$, $P_1$, $P_2$ its three
type~$\sA_6$ singular points, ordered as explained above. Let
$Q_0$, $Q_1$, $Q_2$ be the conics as in Proposition~\ref{prop.conics},
and let $c\colon\Cp2\to\Cp2$ be the
order~$3$ automorphism given by Theorem~\ref{th.symmetry}.
In this section, we assume that $Q_0$, $Q_1$, $Q_2$ are
irreducible, see Lemma~\ref{irreducible}.

Perform the triangular transformation centered at
$P_0$, $P_1$, $P_2$, \ie, blow up the three points and blow down
the proper transforms of the lines $(P_iP_j)$,
$0\le i<j\le2$.
Denote by bars the (proper) images of the curves and points
involved. The curve~$\BB$ has three type~$\sA_4$ singular
points~$\PP_i$, $i=0,1,2$. The transforms~$\QQ_i$, $i=0,1,2$,
are lines, so
that each~$\QQ_i$ passes through~$\PP_{i+1}$ and is tangent
to~$\BB$ at~$\PP_i$. Besides, $\BB$ has three (at least)
nodes~$\SS_i$, $i=0,1,2$, located at the blow-up centers of the
inverse triangular transformation. Note that (under appropriate
indexing) the line $(\SS_i\SS_{i+1})$ contains~$\PP_i$, $i=0,1,2$.

Since the blow-up centers form a single orbit of~$c$, the
transformation commutes with~$c$ and the new configuration is
still $\CG3$-symmetric.

Choose homogeneous coordinates $(u_0:u_1:u_2)$ in~$\Cp2$ so that
$\PP_0=(1:0:0)$, $\PP_1=(0:1:0)$, $\PP_2=(0:0:1)$, and $(1:1:1)$
is one of the fixed points of~$c$. Then $c$
acts \via\ a cyclic permutation of the coordinates, and its
three fixed points are $(1:\Ge:\Ge^2)$, $\Ge^3=1$. The condition that
$\PP_i\in(\SS_i\SS_{i+1})$ and that the triple $\SS_0$, $\SS_1$,
$\SS_2$ is $c$-invariant translates as follows:
there is a parameter~$t\in\C$ such
that $\SS_0=(1:t:t^2)$, $\SS_1=(t^2:1:t)$, and $\SS_2=(t:t^2:1)$.
In order to get three distinct points other than
$\PP_0$, $\PP_1$, $\PP_2$,
one must have $t\ne0$, $t^3\ne1$.

In the chosen coordinates, $\BB$ has three type~$\sA_4$ singular
points located at
the vertices of the coordinate triangle and tangent to its edges.
Since $\BB$ is also preserved by~$c$, it must be given by a polynomial
$F(u_0,u_1,u_2)$ of the form
\begin{equation}\label{eq.abcd}
\deMaple\obeylines\def
{\\}%
\gathered%
a*(u[0]^4*u[2]^2+u[1]^4*u[0]^2+u[2]^4*u[1]^2)+b*(u[0]^3*u[1]^2*u[2]+u[1]^3*u[2]^2*u[0]+u[2]^3*u[0]^2*u[1])
+c*(u[0]^3*u[1]*u[2]^2+u[1]^3*u[2]*u[0]^2+u[2]^3*u[0]*u[1]^2)+d*u[0]^2*u[1]^2*u[2]^2 \endgathered%
\end{equation}
for some $a,b,c,d\in\C$.
Conversely, any curve~$\BB$ given by a polynomial as above
is preserved by~$c$ and has three singular points
adjacent to~$\sA_3$ and situated in the prescribed way with
respect to the coordinate lines~$\QQ_i$, $i=0,1,2$. Due to the
symmetry, it suffices to make sure that $\BB$ is singular
at~$\SS_0$ and that its singularity at~$\PP_0$ is adjacent
to~$\sA_4$. The former condition results in the linear system
\begin{gather*}\displaybreak[0]
6at^4+(3t^4+3t^7)b+(5t^5+t^8)c+2dt^6=0,\\\displaybreak[0]
(4t^3+2t^9)a+(2t^3+4t^6)b+(4t^4+2t^7)c+2dt^5=0,\\
(4t^8+2t^2)a+(t^2+5t^5)b+(3t^3+3t^6)c+2dt^4=0,
\end{gather*}
and the latter condition is equivalent to $b=2a$ or $b=-2a$.
In both cases,
the solution space of the linear system has dimension one; it is spanned by
$$
(a,b,c,d) = (t^2, 2t^2, -2t(t^3+2), (t^3+2)^2)
$$
and
$$
(a,b,c,d) = (1, -2, -2t^2, t(t^3+8)),
$$
respectively.
The first solution, with $b=2a$, results in a reducible polynomial
$$
\deMaple
(t*u[1]^2*u[0]-2*u[1]*u[2]*u[0]-u[1]*u[2]*u[0]*t^3+t*u[0]^2*u[2]+t*u[1]*u[2]^2)^2;
$$
hence, it should be disregarded. For the second solution,
the
substitution $u_0=1$, $u_1=x$, $u_2=y+x^2$ results in a polynomial
with the principal part $y^2-4t^2x^5$, \ie, the singularity
of~$\BB$
at~$\PP_0$ (and, due to the symmetry, at~$\PP_1$ and~$\PP_2$
as well)
is of type~$\sA_4$ exactly. In particular, the curve is
irreducible. (Indeed, the only possible splitting would be into an
irreducible
quintic and a line, but in this case all nodes of~$\BB$ would have
to be collinear.)

To obtain the original curve~$B$, one should perform the
substitution
$$
\deMaple
u[0] = v[0]+t^2*v[1]+t*v[2],\quad
u[1] = t*v[0]+v[1]+t^2*v[2],\quad
u[2] = t^2*v[0]+t*v[1]+v[2]
$$
(passing to an invariant coordinate triangle with the vertices
at~$\SS_i$, $i=0,1,2$) followed by the inverse triangular
transformation $v_0=z_1z_2$, $v_1=z_2z_0$, $v_2=z_0z_1$. The
resulting polynomial is the one given by~\eqref{equation}
with $t\ne0$.

Counting the genus and taking into account the symmetry,
one concludes that the singularities of~$\BB$ at $\SS_i$,
$i=0,1,2$, are either all nodes or all cusps, the latter
possibility corresponding to $t^3=-3$ or
$t^3=-1/3$. Note that the cusps of~$\BB$ merely mean that the original
curve~$B$ is tangent to the lines $(P_iP_j)$, $0\le i<j\le2$;
these curves are still in the same equisingular deformation family.

\begin{remark}\label{triple.conic}
If $t^3=1$, the polynomial~\eqref{equation} becomes reducible. For
example, if $t=1$, it turns into
$\deMaple 4*(z[0]*z[1]+z[1]*z[2]+z[2]*z[0])^3$.
\end{remark}

\subsection{Theorem~\ref{th.equation}: the case of reducible
conics}\label{s.reducible}
Now, assume that the conics $Q_0$, $Q_1$, $Q_2$ are reducible, see
Lemma~\ref{irreducible}. In this case, we can start directly
from~\eqref{eq.abcd}, placing the singular points so that
$P_0=(1:0:0)$, $P_1=(0:0:1)$, and $P_2=(0:1:0)$. Note that
$a\ne0$; hence, we can let $a=1$.

As above, in view of the symmetry it suffices to analyze the
singularity at~$P_0$.

The condition that the singularity is adjacent to~$\sA_4$ is
equivalent to $b=\pm2$. If $b=2$, the substitution
$u_0=1$, $u_1=x$, $u_2=y-x^2+cx^3\!/2$ produces a polynomial
in~$(x,y)$ with the principal part $y^2+(d-c^2\!/4)x^6$. Hence,
$d=c^2\!/4$. However, in this case the original polynomial~$F$ is
reducible:
$$
\deMaple
F=\frac14(2*u[1]^2*u[0]+2*u[2]^2*u[1]+2*u[0]^2*u[2]+c*u[0]*u[1]*u[2])^2.
$$

Let $b=-2$. Then, substituting $u_0=1$, $u_1=x$, $u_2=y+x^2$, one
obtains
$$
y^2+3cyx^3+2cx^5+dx^6-4x^7+\text{(higher order terms)}.
$$
Hence, $c=d=0$, and in this case the singularity at the origin is
exactly~$\sA_6$. The curve is irreducible (as any sextic with
three type~$\sA_6$ singular points) and, after the coordinate
change $(u_0:u_1:u_2)\mapsto(z_0:z_2:z_1)$, the resulting equation
is~\eqref{equation} with $t=0$.

\subsection{Extra singular \term-sextics}\label{s.extra}
Since the total Milnor number of a plane sextic does not
exceed~$19$, a
curve $\curve=\curve(t)$
can have
at most one extra singular point, which must be of type~$\sA_1$.
Since, in addition, $B$ is preserved by~$c$,
this extra singular point must be fixed by~$c$, \ie,
it must be of
the form $(1:\Ge:\Ge^2)$, $\Ge^3=1$. Solving the
corresponding linear system shows that $B$ is singular at
$(1:\Ge:\Ge^2)$, $\Ge^3=1$,
whenever it passes through this point, and this is
the case when $t=-3/\Ge$. In conclusion, $B$ has an extra node
(the set of singularities $3\sA_6\splus\sA_1$)
if and only if $t^3=-27$.

\subsection{Proof of Theorem~\ref{th.equation}}\label{proof.equation}
As shown in Sections~\ref{s.generic} and~\ref{s.reducible}, any
\term-sextic belongs to the connected family~$\curve(t)$,
$t^3\ne1$. According to
Section~\ref{s.extra}, this family represents two equisingular
deformation classes: the restriction to the
connected subset $t^3\ne1,-27$ (the set of singularities $3\sA_6$)
and three equivalent isolated curves corresponding to
$t^3=-27$ (the set of singularities $3\sA_6\splus\sA_1$;
the
equivalence is given by the
coordinate
change
$(z_0:z_1:z_2)\mapsto(z_0:\Ge z_1:\Ge^2z_2)$, $\Ge^3=1$, \cf.
Remark~\ref{rem.triple}). Comparing this result with
Proposition~\ref{prop.sing}, one concludes that any curve given
by~\eqref{equation} with $t^3\ne1$ is a \term-sextic.
\qed

\section{The fundamental group}\label{S.group}

\subsection{Calculation of the
group}

For the calculation, we choose a real curve~$\curve(t)$ given by
Theorem~\ref{th.equation} and close to the triple conic~$\curve(1)$,
see Remark~\ref{triple.conic}.

\begin{theorem}\label{th.epi}
For the curve $\curve=\curve(5/6)$, there is an epimorphism
$$
G:=\langle\Go,\xi\,|\,\xi^2=\1,\ \Go^2=\xi\Go^5\xi\rangle
 \twoheadrightarrow\pi_1(\Cp2\sminus\curve).
$$
\end{theorem}

\begin{proof}
Let $\curve=\curve(5/6)$ and make the following change of coordinates:
\[
z_0= 1/3\,{ Z}-1/3\,Y+1/3\,X,\,\,z_1=-1/3\,X+2/3\,{ Z}-5/3\,Y,\,\,
 z_2=Y.
\]
Then, in the affine space $\C^2=\Cp2\sminus\{Z=0\}$, $\curve$ is defined by
$g(x,y)=0$, where
\begin{multline*}
g(x,y)={\frac {716}{19683}}\,x+{\frac {17872}{177147}}\,y-{\frac {11503}{
708588}}\,xy-{\frac {356093}{354294}}\,x{y}^{2}+{\frac {322559}{
5668704}}\,{x}^{3}y\\
+{\frac {3568}{177147}}-{\frac {722513}{2834352}}\,
{x}^{2}y-{\frac {8137}{472392}}\,{x}^{2}-{\frac {28582655}{2834352}}\,
x{y}^{5}+{\frac {56261293}{22674816}}\,{y}^{4}\\
-{\frac {449027}{1417176
}}\,{y}^{2}-{\frac {4427549}{2834352}}\,{y}^{3}+{\frac {81485377}{
11337408}}\,{y}^{5}-{\frac {255219619}{22674816}}\,{y}^{6}\\
-{\frac {
57539}{1417176}}\,{x}^{3}+{\frac {2243}{209952}}\,{x}^{4}-{\frac {
26011}{5668704}}\,{x}^{6}+{\frac {2726579}{7558272}}\,{x}^{2}{y}^{2}\\
+{
\frac {1092623}{1259712}}\,x{y}^{3}+{\frac {9868757}{11337408}}\,{x}^{
3}{y}^{2}+{\frac {11718893}{3779136}}\,{x}^{2}{y}^{3}+{\frac {77768419
}{11337408}}\,x{y}^{4}\\
-{\frac {309307}{5668704}}\,y{x}^{5}-{\frac {
9030539}{22674816}}\,{x}^{4}{y}^{2}-{\frac {9923629}{5668704}}\,{x}^{3
}{y}^{3}\\
-{\frac {61362175}{11337408}}\,{x}^{2}{y}^{4}+{\frac {12505}{
944784}}\,{x}^{5}+{\frac {397175}{2834352}}\,{x}^{4}y.
\end{multline*}
Its three singularities are located at
\[
 P_0=(-1,0),\quad P_1=(2,0),\quad P_2=(-1/2,1/2),
\]
and the graph in the real plane is given
in
Figure \ref{graph}.

\begin{figure}[htb,here]
\setlength{\unitlength}{1bp}
\begin{picture}(600,250)(-100,0)
\scalebox{0.6}[0.4]
{\includegraphics{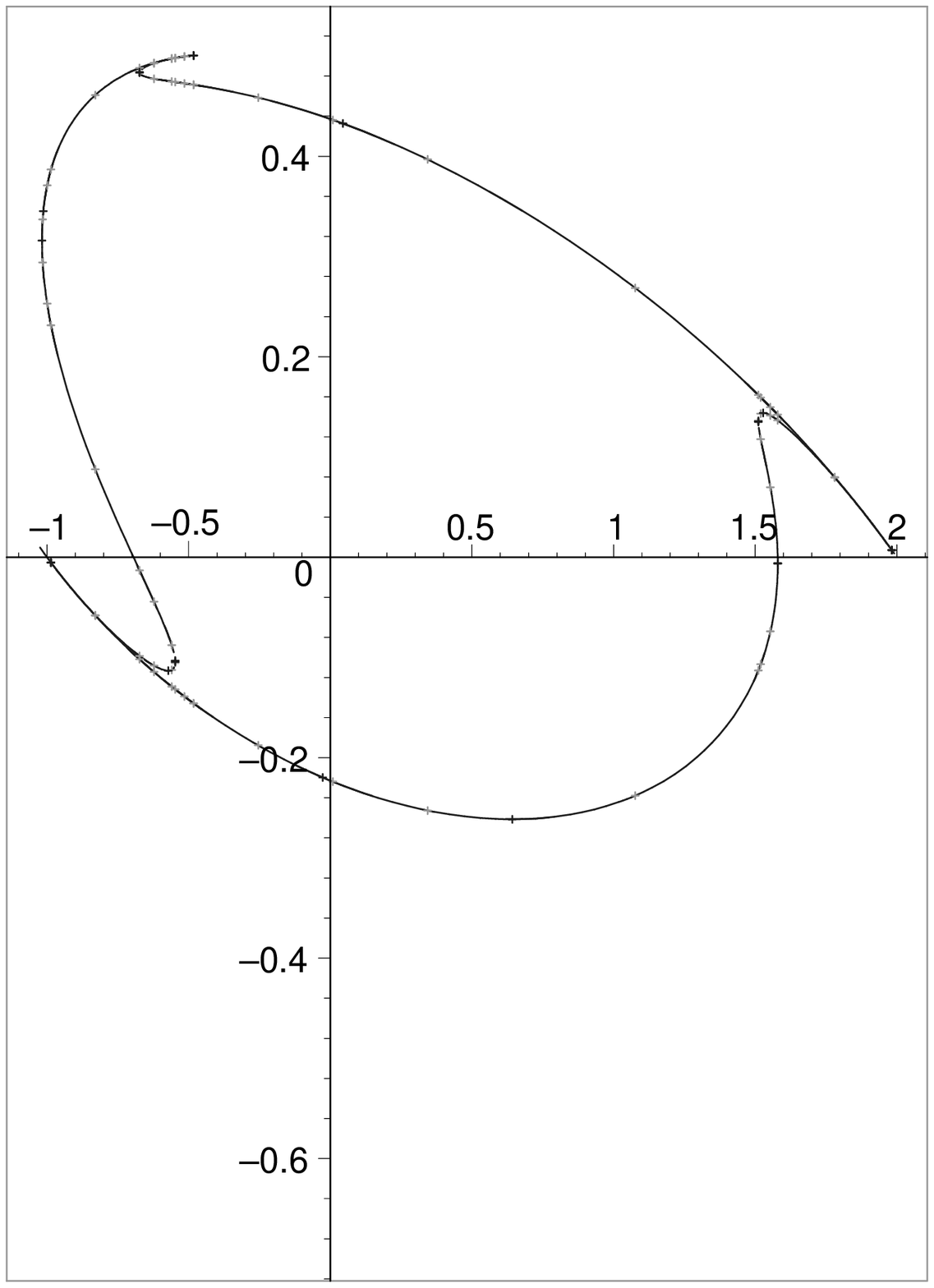}}
\end{picture}
\caption{The graph of $\curve$}\label{graph}
\end{figure}

For the calculation, we
apply van Kampen's method~\cite{vanKampen} to the
horizontal pencil
$L_\eta=\{y=\eta\}$. The singular pencils corresponds to the roots of
\begin{multline*}
\!\!( 90617210907008\,{y}^{9}-60741238168704\,{y}^{8}-
52338630572904\,{y}^{7}+38781803208839\,{y}^{6}\\
+8841431367018\,{y}^{5}
-8143800845364\,{y}^{4}-176669916264\,{y}^{3}+512733413664\,{y}^{2}\\-
7789219200\,y-6298560000 )
{y}^{14}  \left( 2\,y-1 \right) ^{7}=0.
\end{multline*}
Note that we have five real singular pencil lines
$$
\gathered
 L_\eta,\quad \eta=\eta_i,\,i=1,2,\dots, 5,\\
\eta_1\approx -0.26,\quad\eta_2\approx -0.11,\quad
 \eta_3=0,\quad\eta_4\approx 0.14,\quad \eta_5=1/2,
\endgathered
$$
where $L_{\eta_i},\,i=1,2,4$ are tangent to $\curve$, which come from
 the first factor of degree 9.
There also are three pairs of complex conjugate singular fibers,
but we do not use them;
that is why we only assert that the
map constructed below is an epimorphism, not an isomorphism.
We  take the base point at infinity $b=(1:0:0)$,
and we fix generators $\rho_1,\dots,\rho_6$
on the regular pencil line
$y=-\eps$ (where $\eps$ is a sufficiently small positive real number)
as in Figure~\ref{Generators}. The bullets are lassos, which
are counterclockwise oriented loops
going around a point of~$\curve$.

\begin{figure}[htb]
\setlength{\unitlength}{1bp}
\begin{picture}(600,150)(-100,0)
\includegraphics[width=6cm,clip]{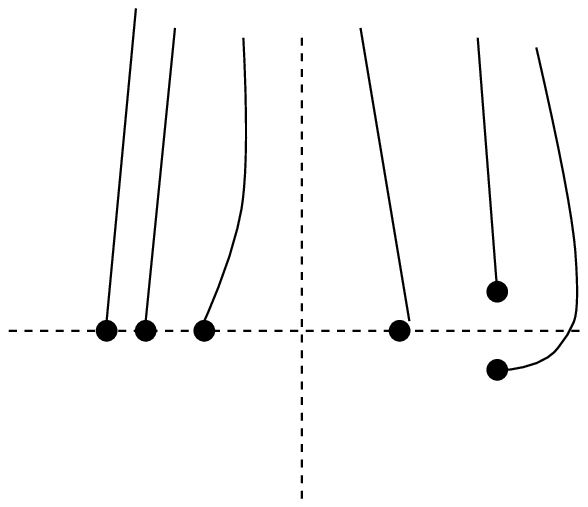}
\put(-150,120){$\rho_6$}
\put(-120,120){$\rho_5$}
\put(-80,130){$\rho_3$}
\put (-100,120){$\rho_4$}
\put(-50,135){$\rho_2$}
\put(-10,125){$\rho_1$}
\put(-80,45){$O$}
\end{picture}
\caption{Generators in the fiber $y=-\eps$}\label{Generators}
\end{figure}

The monodromy relations at $y=\eta_2$, $\eta_1$ are
tangent relations,
they are given as
\rel(R_1):\rho_4=\rho_5,\quad \rho_3=\rho_4^{-1}\rho_6\rho_4\endrel
Thus, hereafter we eliminate the generator~$\rho_5$.

The monodromy relations at $y=0$ are two $\sA_6$-cusp relations,
they are given as
\rel(R_2):\begin{cases}
\omega^3\rho_6=\rho_4\omega^3,\quad& \omega=\rho_6\rho_5,\\
\tau^3\rho_2=\rho_1\tau^3,\quad & \tau=\rho_2\rho_1.
\end{cases}
\endrel
To see the relations at $y=\eta_4$, $\eta_5$ effectively, we take
new elements $\rho_6'$, $\rho_5'$, $\rho_2'$, $\rho_1'$ as in
Figure~\ref{y=eps}.
The
new elements are defined as
\rel(R_3):
 \rho_6'=\omega^{-2}\rho_6\omega^2,\quad
 \rho_4'=\omega^{-1}\rho_4\omega,\quad
 \rho_2'=\tau^{-1}\rho_1\tau,\quad
 \rho_1'=\tau^{-2}\rho_2\tau^2.
\endrel
Note that
they satisfy the relations
\[
 \rho_4'\rho_6'=\omega,\quad \rho_2'\rho_1'=\tau.
\]

\begin{figure}[htb]
\setlength{\unitlength}{1bp}
\begin{picture}(600,150)(-100,0)
\includegraphics[width=7cm,clip]{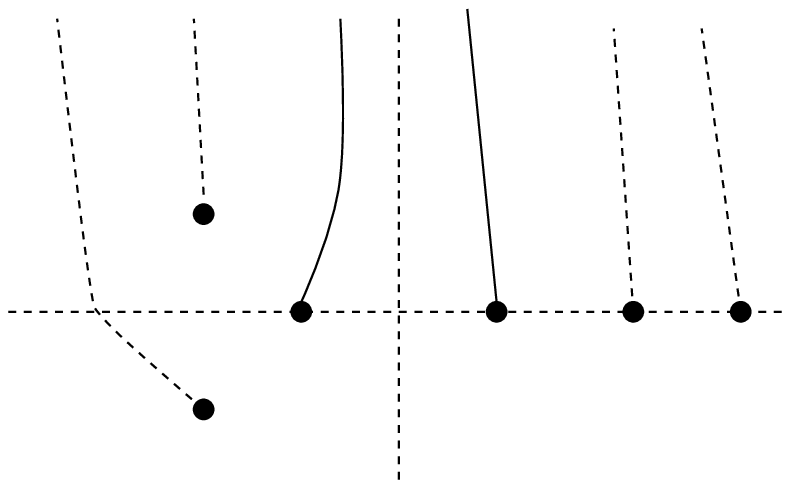}
\put(-180,110){$\rho_4'$}
\put(-150,110){$\rho_6'$}
\put (-125,100){$\rho_4$}
\put(-95,105){$\rho_3$}
\put(-55,105){$\rho_2'$}
\put(-35,105){$\rho_1'$}
\end{picture}
\caption{Generators in the fiber $y=\eps$}\label{y=eps}
\end{figure}

Now, the monodromy relation at $y=\eta_4$
is given as
\rel(R_4):
\rho_3=\rho_2'\quad \text{or}\quad
\rho_3=\tau^{-1}\rho_1\tau.
\endrel
The relation at $y=1/2$ is an $\sA_6$-cusp relation,
which is given as
\rel(R_5):
(\rho_4\rho_2'\rho_1'{\rho_2'}^{-1})^3\rho_4=(\rho_2'\rho_1'{\rho_2'}^{-1})
(\rho_4\rho_2'\rho_1'{\rho_2'}^{-1})^3.
\endrel
Finally, the vanishing relation at infinity is given as
\rel(R_{\infty}): \omega^2\tau=\1.\endrel

We eliminate the generator~$\rho_3$
using~$(R_1)$.
Then $(R_4)$ is translated into the following
relation:
\rel(R_4'): \rho_4^{-1}\rho_5\rho_4=\tau^{-1}\rho_1\tau.\endrel
The relation $(R_\infty)$ can be rewritten as
\rel(R_\infty'): \tau=\omega^{-2}.\endrel
From $(R_\infty')$ and $(R_4')$, we get
\[
 \rho_1=(\tau\rho_4^{-1})\rho_6(\rho_4\tau^{-1})
=\omega^{-2}\rho_4^{-1}\rho_6\rho_4\omega^2
\underset{R_2}=\rho_6\rho_4\rho_6^{-1}.
\]
As $\rho_2=\tau\rho_1^{-1}$, this implies
\rel(R_4''):\rho_1=\rho_6\rho_4\rho_6^{-1},\quad
\rho_2=
\omega^{-1}\rho_4^{-1}\omega^{-1}.
\endrel
We can rewrite $\rho_2'$, using the above relations, as
follows:
\[
 \rho_2'=\omega^3\rho_6^{-1}\omega^{-2}\underset{R_2}= \rho_4^{-1}\rho_6\rho_4.
\]
Thus, $\rho_2'\rho_1'{\rho_2'}^{-1}=\omega^{-3}\rho_4$, and $(R_5)$
can be rewritten in $\rho_6$, $\rho_5$ as follows:
\rel(R_5'):
(\rho_4\rho_6\omega^{-3})^3\rho_4=(\rho_6\omega^{-3})(\rho_4\rho_6\omega^{-3})^3.
\endrel
We have to rewrite the relations in the words of $\rho_4,\rho_6$.
The relation $\tau^3\rho_2=\rho_1\tau^3$ gives
\[
 \omega^{-6}(\omega^{-1}\rho_4^{-1}\omega^{-1})=(\rho_6\rho_4\rho_6^{-1})\omega^{-2},
\]
which reduces to $\omega^5\rho_6=\rho_4\omega^8$.
Using the relation $\omega^3\rho_6=\rho_4\omega^3$ several times, we get
$\rho_4=\omega\rho_6\omega^2$.
Now, we eliminate $\rho_4$ using $\omega=\rho_6\rho_4$ to obtain
\rel(R_6):(\omega\rho_6)^2=\1.\endrel
Replace
the generator~$\rho_6$ by $\xi=\omega\rho_6$,
so that the new generators are $\omega$, $\xi$;
then $\rho_4$, $\rho_6$ are
expressed as
  $\rho_6=\omega^{-1}\xi$ and $\rho_4=\xi^{-1}\omega^2$,
and $(R_6)$ is written as $\xi^2=\1$.
The relation $\omega^3\rho_6=\rho_4\omega^3$ reduces to
\rel(R_7): \omega^2=\xi\omega^5\xi.\endrel
Thus, we have shown that $\pi_1(\Cp2\sminus\curve)$
is generated by two elements $\omega$, $\xi$,
which are subject to the relations
\[
 \xi^2=\1,\quad \omega^2=\xi\omega^5\xi.
\]
This establishes the required epimorphism.
\end{proof}

\subsection{The group structure of $G$}
Below,
we analyze the group~$G$ obtained in Theorem~\ref{th.epi} and show
that it is isomorphic to
\[
 \DG{14}\times\CG3=
 \langle a,b,\xi\,|\,\xi^2=\1,\ a^3=\1,\ b^7=\1,\ \xi b\xi=b^6,\
 [a,b]=[a,\xi]=\1
\rangle,
\]
where $[a,b]=aba^{-1}b^{-1}$ is the commutator.

\begin{lemma}\label{G=DZ}
The map $\xi\mapsto\xi$, $\omega\mapsto ab$ establishes an
isomorphism $G\cong\DG{14}\times\CG3$.
\end{lemma}

\begin{proof}
Putting $c=ab$, we see that
$a=c^7,\,b=c^{15}$. Thus we can use two generators $\xi$, $c$, and
\[
  \DG{14}\times\CG3=\langle \xi,c\,|\,
\xi^2=\1,\ c^{21}=\1,\ \xi c^{15}\xi=c^6,\ [c^7,\xi]=\1\rangle.
\]
Now we consider our group~$G$:
\[
 G=\langle \omega,\,\xi\,|\,\xi^2=\1,\ \xi\omega^5\xi=\omega^2
\rangle.
\]
First we see that
$\xi\,\omega^{25}\,\xi=\omega^{10}=(\xi\,\omega^2\,\xi)^2=\xi\omega^{4}\xi$.
Thus we get $\omega^{21}=\1$.
We assert that $\omega^7$ is in the center of $G$.
Indeed,
$\xi\omega^{14}\xi=(\xi\,\omega^2\,\xi)^7=\omega^{35}=\omega^{14}$.
Thus, $\omega^{14}$ is in the center,
and so is
$\omega^7=(\omega^{14})^2$.
Observe that
$\xi\omega^{15}\xi\,=\, (\omega^2)^3\,=\,\omega^6$.
Thus, we have another presentation of $G$,
\[
 G=\langle\,\omega,\,\xi\,|\, \xi^2=\1,\
\omega^{21}=\1,\ \xi\omega^{15}\xi=\omega^6,\ [\omega^7,\xi]=\1
\rangle,
\]
which coincides with that of $\DG{14}\times\CG3$.
(The original relation $\xi\omega^2\xi=\omega^5$ is recovered
by squaring the relation
$\xi\omega^{15}\xi=\omega^6$, taking
into account $\Go^{21}=\1$, and
cancelling the central element~$\Go^7$.)
\end{proof}


\subsection{Proof of Theorem~\ref{th.main}}\label{proof.main}
According to Theorem~\ref{th.equation},
any curve $\curve=\curve(t)$,
$t^3\ne1$, is a \term-sextic, \ie, its fundamental group
$\pi=\pi_1(\Cp2\sminus\curve)$ factors to~$\DG{14}$. On the other
hand,
$\pi/[\pi,\pi]=\CG6$. The smallest group with these properties
is $\DG{14}\times\CG3$, \ie,
one has $\ord\pi\ge\ord(\DG{14}\times\CG3)$.
In view of Theorem~\ref{th.epi} and Lemma~\ref{G=DZ}, there is
an epimorphism $\DG{14}\times\CG3\twoheadrightarrow\pi$; comparing
the orders, one concludes that it
is an isomorphism.
\qed

\end{document}